\RequirePackage{fix-cm}

\documentclass[numbook]{svjour3}

\usepackage{amsmath,amsxtra,amssymb, amscd,amsfonts,enumerate}
\usepackage[mathscr]{eucal}
\usepackage{graphicx} 
\usepackage[misc]{ifsym}

\usepackage{graphicx,color}

\DeclareMathOperator{\R}{\mathbb{R}}

\allowdisplaybreaks
\journalname{Journal Name}

\begin{document}

\titlerunning{Viscosity solutions of Hamilton-Jacobi equations in Banach spaces}

\title{On the existence, uniqueness and stability of  $\beta$-viscosity solutions to a class of Hamilton-Jacobi equations in Banach spaces}

\author{Tran Van Bang \and Phan Trong Tien} 

\institute{
\begin{itemize}
\item[] 2010 \textit{Mathematics Subject Classification}. 35R15, 49J52, 35D40.
\item[]This research was supported by the Ministry of Education and Training of Vietnam under Grant No. B2018-SP2-14.\\ 
\item[\Letter] T.V. Bang\\
tranvanbang@hpu2.edu.vn\medskip
\item[] Department of Mathematics, Hanoi Pedagogical University 2, Vinh Phuc, Viet Nam
\item[\Letter] P.T. Tien\\
trongtien2000@gmail.com\medskip
\item[] Department of Mathematics, Quang Binh University, Quang Binh, Viet Nam
\end{itemize}
}

\date{Received: 25 February 2018 / Accepted for publication in AMV: 06 June 2018}

\maketitle

\begin{abstract}
This paper is concerned with the qualitative properties of viscocity solutions
to a class of Hamilton-Jacobi equations (HJEs) in Banach spaces.
Specifically, based on the concept of $\beta$-derivative \cite{DGZ93b}
we establish the existence, uniqueness and stability of  $\beta$-viscosity solutions
for a class of HJEs in the form $u+H(x,u,Du)=0.$
The obtained results in this paper extend ealier works
in the literature, for example, \cite{CL85}, \cite{CL86} and \cite{DGZ93b}.
\keywords{$\beta$-borno, viscosity solutions, Hamilton-Jacobi equations, nonlinear partial differential equations.} 
\end{abstract}

\section{Introduction}
Let $X$ be a real Banach space with a $\beta$-smooth norm $|\cdot|$ (see Section 2 for more details), $\Omega\subset X$ an open subset. We study the existence, uniqueness and stability of $\beta$-viscosity solutions for the following HJEs
\begin{equation}\label{A1}
	u+ H(x,u,Du)=0\text{ in }\Omega,
\end{equation}  
subject to the boundary condition (in the case $\Omega\not=X$)
 \begin{equation}\label{A2}
 u=\varphi \text{ on }\partial \Omega.
 \end{equation}
Here, $\varphi:\partial\Omega\to \R$ and $H:\overline{\Omega}\times \R\times X^*_\beta\to \R$ are merely continuous in general, where $X^*_\beta$ is the dual space of the Banach space $X,$ and equipped with topology $\tau_{\beta}$ (see Definition 2.2). 
   
Hamilton-Jacobi equations (HJEs) comprise an important class of
nonlinear partial differential equations (PDEs). This type of PDEs
are used to describe dynamics of various practical models and physical processes.
Typical examples of using HJEs in such models can be found
in mechanics or optimal control theory, especially for PDEs.
For instance, in the optimal control theory, a Bellman equation,
also known as a dynamic programming equation, is a necessary condition
for optimality associated with the dynamic programming,
which can be characterized by the so-called Hamilton-Jacobi-Bellman equation,
a subclass of HJEs.  These equations in general, are the nonlinear ones and  often do not have classical solutions, so it is necessary to study weak solutions in some sense.

The concept of viscosity solutions of nonlinear PDEs
was first introduced for HJEs in finite-dimensional spaces
in the earlier of 1980s \cite{dau}. Due to its adaptive capability
in practical applications, particularly in optimal control theory,
this concept has received considerable research attention.
We refer the reader to \cite{AR,H2013} and the references therein.
Recently, the studies on viscosity solutions have also
been developed for HJEs in infinite-dimensional spaces \cite{BZ99,BZ96,CL86,CL85,DGZ93,DGZ93b}
and second-order PDEs both in finite-dimensional spaces \cite{CIL92}
and infinite-dimensional spaces \cite{B2006,BV2006}.
It should be noted that viscosity solutions of PDEs
are weak solutions since they are merely continuous functions whose derivatives
are defined through test functions using extreme principles.
However, it has been shown that some specific viscosity tests
can be conducted through upper and lower differentiations, which are also known as semi-derivatives. 
This creates a close connection between the theory of viscosity solutions
and non-smooth analysis including the theory of subdifferentials.

It is known that numerous types of subderivatives (superderivatives)
have been proposed in the literature. Among that the subderivatives
in the sense of  Fr\'{e}chet, Hadamard, G\^{a}teaux and Mordukhovich
are the most widely used ones \cite{BT2014,BZ99,DGZ93b,D2003,MNY2007,MYZ2000}.
Clearly, for a class of HJEs, the use of different subderivatives
may lead to different types of viscosity solutions.
In pioneering works  \cite{AR,CL86,CL85,dau}, the meaning of viscosity solutions
is characterized by semi-Fr\'{e}chet derivative.
On the other hand, in many existing works,
the study on qualitative properties of viscosity solutions is
usally based on a certain type of smoothness of the underlying norm.
However, such smoothness is normally not available in most popular Banach spaces like $L^1.$
To overcome this issue, the authors of \cite{BZ96,DGZ93b}
proposed the concepts of Borno $\beta,$ derivative $\beta$-viscosity, $\beta$-viscosity solutions and obtained a unique result for the HJEs of the form $u+H(x,Du)=0.$

In this paper, we extend the unique result of \cite{BZ96} to a broader class of equations \eqref{A1}, and also  establish the existence and stability of $\beta$-viscosity solutions under certain assumptions. In the literature, 
the smooth variational principle proved by Deville in [12] has been employed as a key tool in proving the uniqueness of $\beta$-viscosity solutions of Hamilton-Jacobi equations of the form $u+F(Du)=f$, where $F$ is uniformly continuous on $X^*_\beta$ and $f$ is uniformly continuous and bounded on $X$. However, it is nontrivial to apply that principle for equations of the general form $u+H(x, u, Du)=0$ on a set $\Omega\subset X$. In [8], Crandall and Lions  
established the uniqueness of Frechet viscosity solutions for Hamilton-Jacobi equations of the above general form by using the Radon-Nikodym property as a key assumption. By the Radon-Nikodym property, we means that if $\varphi$ is a bounded and lower semicontinuous real-valued function on a closed
ball $B$ in $X$ and $\varepsilon > 0$, then there is an element $x^*\in X^*$
of norm at most $\varepsilon$ such
that $\varphi+x^*$
attains its minimum on $B$. In this work, we use certain extensions of functions to prove the smooth variational principle for functions defined on a subset $E$ of $X$. Then we will establish the uniqueness of $\beta$-viscosity solutions (which are weaker than Fr\'{e}chet-viscosity solutions)   for $u+H(x, u, Du)=0$ on a set $\Omega\subset X$  by the doubling of variables technique. Moreover, we work on a Banach space $X$ having a $\beta$-smooth norm or a norm being equivalent to a $\beta$-smooth norm. The Radon-Nikodym assumption is unnecessary in our arguments. Another contribution of our paper is on the existence and stability of solutions.  More specifically, by Perron's method, we show that under certain conditions, the equations of the form $u+H(x, u, Du)=0$ have $\beta$-viscosity solutions and these solutions as well as the stability of this class of equations.

The remaining of this paper is organized as follows.
Some preliminaries on $\beta$-viscosity subderivatives,
smooth variational principle on Banach spaces and discussion about our main assumptions on the ambient space are presented in Section 2.
Our main results are presented in Section 3. The key result is stated in Theorem \ref{dl3.1} whose proof is postponed to Section 4. From this, the comparison and uniqueness of solutions are proved, as a consequence (Corollary \ref{hq3.1}). The obtained results extend those of  \cite{CL85}, to spaces
endowed with $\beta$-smooth norms, equivalent $\beta$-smooth norms
or the ones satisfying the hypothesis  $(H_{\beta})$ without Radon-Nikodym property. The existence and stability are respectively showed in Theorem \ref{dl3.2} and Theorem \ref{dl3.3}. These are the extension of those in  \cite{CL86} to spaces
endowed with $\beta$-smooth norms instead of Fr\'echet-borno. As mentioned above, Section 4 is used for the proof of Theorem \ref{dl3.1}. In the last section, we give the conclusion and some respects for future study. 
\section{Preliminaries}
Let $X$ be a real Banach space with the norm $|.|.$   
The close unit ball and the dual space of $X$ will be
denoted as $B_1$ and $X^*$, respectively. We denote by $\left< x^*,x \right>$
the value of  $x^*\in X^*$ at $x\in X;$ that is,  $\left< x^*,x\right> =x^*(x).$ The norm on the dual space $X^*$ will be also denoted as $|.|.$   

Throughout this paper, the concerned functions are defined on
the Banach space $X$ taking values in the extended real line $\overline{\mathbb{R}}.$
In addition, it is our convention here that when a function
is lower semicontinuous (resp. upper semicontinuous), it is not identical to $-\infty$ (resp. $+\infty$).

\begin{definition}
A \emph{borno} $\beta$ on $X$ is  a family of closed, bounded, and centrally symmetric subsets of $X$
satisfying the following three conditions:

a) $X=\bigcup\limits_{B\in \beta}B,$

b) $\beta$ is closed under scalar multiplication,

c) the union of any two elements in $\beta$  is contained in some element of $\beta.$
\end{definition}

\begin{example}\label{vd01}
It is easy to verify the following facts.
\begin{enumerate}[1)]
\item The family $F$ of all closed, bounded, and centrally symmetric subsets of $X$ is a borno
on $X,$ which is called \textit{Fr\'{e}chet borno}.
\item The family $H$ of all compact, centrally symmetric subsets of $X$
is a borno on $X$ called \textit{Hadamard borno}.
\item The family $WH$ of all weakly compact, closed, and centrally symmetric subsets of $X$
is a borno on $X$ called \textit{weak Hadamard borno}.
\item The family $G$ of all finite, centrally symmetric subsets of $X$ is also a borno on $X$ called \textit{G\^{a}teaux borno}.
\end{enumerate}
\end{example}

\begin{definition}
A sequence $(f_m)_m\subset X^*$ is said to converge to an $f\in X^*$
in $\beta$-topology $\tau _ {\beta}$ if  $f_m\to f $ when $m\to \infty$
uniformly over all elements of $ \beta,$ this means that
for each set $ M \in \beta $ and any given $\varepsilon>0,$
there exists an $n_0\in \mathbb{N}$ such that $|f_m (x) -f (x)| <\varepsilon$
holds for all $m \ge n_0$ and $x\in M.$
\end{definition}

Suppose that $\beta$ is a  borno  on $X.$ We denote by $\tau_{\beta}$
the topology on $X^*$ associated with the uniform convergence
on $\beta$ subsets and $X^*_{\beta}$ the topological vector space $(X^*, \tau_\beta).$

\begin{remark}\label{nx02}
If  $\beta$ borno  is  $F$  (Fr\'{e}chet), $H$ (Hadamard), $WH$ (Hadamard weak) or $G$ (G\^{a}teaux),
then we have Fr\'{e}chet topology, Hadamard topology, Hadamard weak topology and
G\^{a}teaux topology on the dual space $X^*,$ respectively \cite{DGZ93b}.
Thus, $F$-topology is the strongest topology and $G$ topology is the weakest topology
among $\beta$-topologies on $X^*.$
\end{remark}

\begin{definition}[\cite{BZ96}]
Given a function $f: X \to \overline{\mathbb{R}}.$
We say that $f$ is $\beta$-\textit{differentiable} at $x_0\in X$
with $\beta$-\textit{derivative} $\nabla_\beta f(x_0) = p \in X^*$ if $f(x_0) \in\mathbb{R}$ and
\begin{equation*}\label{eq:01}
\lim_{t\to 0}\frac{f(x_0+th)-f(x_0)-\left<p,th\right>}{t}=0
\end{equation*}
uniformly in $h\in V$ for every $V\in\beta.$

We say that the function $f$ is $\beta$-smooth at $x_0$
if there exists a neighborhood $U$ of $x_0$ such that $f$ is $\beta$-differentiable on $U$
and $\nabla_\beta f: U \to X^*_{\beta}$ is continuous.
\end{definition}

\begin{definition}[\cite{BZ96}]
Let $f: X\to\overline{\R}$ be a lower semicontinuous function and
$f(x)<+\infty.$ We say that $f$ is $\beta$-\textit{viscosity subdifferentiable}
and $x^*$ is a $\beta$-\textit{viscosity subderivative} of $f$ at $x$
if there exists a local Lipschitzian function $g: X\to \mathbb{R}$ such that
$g$ is $\beta$-smooth at $x,$ $ \nabla_\beta g (x) = x^*$ and
$f-g$ attains a local minimum at $ x.$
We denote the set of all $\beta$-subderivatives of $f$ at $x$
by $D^-_{\beta}f(x),$ which is called \textit{$\beta$-viscosity subdifferential} of $f$ at $x.$

Let $f: X\to\overline{\R}$ be an upper semicontinuous function and $f(x)>-\infty.$
We say that $f$ is $\beta$-\textit{viscosity superdifferentiable}
and $x^*$ is a $\beta$-\textit{viscosity superderivative} of $f$ at $x$ if
there exists a local Lipschitzian function $g: X\to \mathbb{R}$ such that
$g$ is $\beta$-smooth at $x,$ $ \nabla_\beta g (x) = x^*$ and $f-g$ attains
a local maximum at $ x.$ We denote the set of all $\beta$-superderivatives of
$f$ at $x$ by $D^+_{\beta}f(x),$ which is called \textit{$\beta$-viscosity superdifferential} of $f$ at $x.$
\end{definition}

Next, we denote 
$$\mathcal{D}_{\beta}(X)=\{g: X\to \R| \ g \text{ is bounded, Lipschitzian, and }\beta\text{-differentiable on } X\},$$ 
$$\|g\|_{\infty}= \sup\{|g(x)|: \ x\in X\},\quad \|\nabla_\beta g\|_{\infty}=\sup\{\|\nabla_\beta g(x)\|: \ x\in X\}$$
and
$$ \mathcal{D}^*_{\beta}(X)=\{g\in \mathcal{D}_{\beta}(X)| \ \nabla_\beta g:X\to X_{\beta}^* \text{ is  continuous}\}.$$

The following hypotheses will be used in the derivation of our results.
\begin{itemize}
\item[$(H_{\beta})$] There exists a  bump function $b$ such that $b\in \mathcal{D}_{\beta}(X);$ and
\item[$(H^*_{\beta})$] There exists a  bump function $b$ (i.e. its support is nonempty and bounded)
such that $b\in \mathcal{D}^*_{\beta}(X).$
\end{itemize}

\begin{proposition}\label{mdkg}
The hypotheses $(H_{\beta}) $ and $(H^*_{\beta})$ are fulfilled
if the Banach space $X$ has a $\beta$-smooth norm.
\end{proposition}

\begin{proof}
Let $f: \mathbb{R} \to \mathbb{R}$ be a continuously differentiable function
with bounded derivative whose support is confined in the interval $[1,3].$
For instance, such a function $f$ is given by
$$f(t)=\begin{cases}\exp\left(\frac{1}{(t-2)^2-1}\right)& \textit{ if }\  |t-2|< 1,\\
0&\textit{ if }\  |t-2|\geq1.
\end{cases}$$ 
We define the function $b(x)=f(|x|)$ for $x\in X.$
Clearly $b$ is bounded and Lipschitz continuous.
For $|x|\leq 1$ or $|x|\geq3,$ $b(x)=0.$ Due to the $\beta$-smoothness of the norm,
$b$ is differentiable at all points satisfying $1<|x|<3.$
Combining with the smoothness of  $f,$ we can conclude that $(H_{\beta}) $ and $(H^*_{\beta})$ are satisfied.\hfill $\square$
\end{proof}

In general, the converse of
Proposition \ref{mdkg} does not hold. A counterexample can be found in \cite[Page 59]{Ha2},
where the author pointed out that
there exists a Lipschitz continuous bump function, which is Fr\'{e}chet-differentiable
 in a Banach space $X$ without any  G\^{a}teaux-differentiable norm (hence, there does not exist Fr\'{e}chet-differentiable norm).

In \cite{DGZ93}, the authors established a smooth variational principle in a Banach space $X$
where there exist Lipschitz, bounded and differentiable bump functions.
That principle is stated in the following proposition.

\begin{proposition}[Remark II.5, \cite{DGZ93}]\label{dlbp}
Let $X$ be a Banach space satisfying the hypotheses $(H_{\beta})$ (resp. $(H^*_{\beta})$),and
$f$ a lower semicontinuous function on $X.$
Then there  exists a constant $a=a_X>0$
such that for any $\varepsilon \in (0,1)$ and $y_0\in X$ satisfying $f(y_0)<\inf_{X}(f)+a\varepsilon^2,$
there exist a $g\in \mathcal{D}_{\beta}(X)$ (resp $g\in\mathcal{D}^*_{\beta}(X)$) and an $x_0\in X$ such that:
\begin{enumerate}[(a)]
	\item $f+g$ attains its minimum at $x_0.$
	\item $\|g\|_{\infty}\leq \varepsilon$ and $\|\nabla_{\beta}g\|_{\infty}\leq \varepsilon.$
	\item $|x_0-y_0|\le \varepsilon.$
\end{enumerate}
\end{proposition}

From Proposition \ref{dlbp} we have the following result.

\begin{proposition}\label{dlbp1}
Let $X$ be a Banach space satisfying $(H_{\beta})$ (resp. $(H^*_{\beta})$)
and $E$ a closed subset of  $X.$ Then, for a lower semicontinuous bounded from below function $f$
on $E$ and any $\varepsilon \in (0,1),$ there exist
a $g\in \mathcal{D}_{\beta}(X)$ (resp. $g\in\mathcal{D}^*_{\beta}(X)$) and an $x_0\in E$ such that:
\begin{enumerate}[(a)]
	\item $f+g$ attains its minimum at $x_0.$
	\item $\|g\|_{\infty}\leq \varepsilon$ and $\|\nabla_{\beta}g\|_{\infty}\leq \varepsilon.$
\end{enumerate}
\end{proposition}

\begin{proof}
We extend the function $f$ to the function $\widetilde{f}$ defined on $X$ by
$$\widetilde{f}(x)=\begin{cases}	f(x)& \text{ if } x\in E,\\	\infty& \text{ if } x\notin E.\end{cases}$$
Since $E$ is  closed, it follows from the properties of $f$ that
$\widetilde{f}$ is lower semicontinuous and bounded from below.
Thus, for any $\varepsilon>0,$ there exists an $y_0\in X$
such that $\widetilde{f}(y_0)<\inf_X(\widetilde{f})+\varepsilon.$
Moreover, $\widetilde{f}(y_0)<\infty$ and $y_0\in E$ since $\inf_X(\widetilde{f})+\varepsilon<\infty.$
According to Proposition \ref{dlbp}, there exist a $g\in \mathcal{D}_{\beta}(X)$
(resp. $g\in\mathcal{D}^*_{\beta}(X)$) and an $x_0\in X$ such that:
\begin{enumerate}[(a)]
\item $\widetilde{f}+g$  attains its minimum at $x_0.$
\item $\|g\|_{\infty}\le \varepsilon$ and $\|\nabla_{\beta}g\|_{\infty}\le \varepsilon.$
\end{enumerate}

It follows from (a) that $(\widetilde{f}+g)(y)\geq (\widetilde{f}+g)(x_0)$
for any $y\in X.$ In particular, for $y\in E,$ we have $\infty> (\widetilde{f}+g)(y)=(f+g)(y)\geq (\widetilde{f}+g)(x_0).$
Since $\|g\|_{\infty}\leq \varepsilon,$ we can conclude that $x_0\in E.$
The proof is completed.\hfill $\square$
\end{proof}

Similarly, when $f$  is upper semicontinuous and bounded from above on $E,$
the function $\widehat{f}=-f$ is lower semicontinuous and bounded from below on $E.$
We have the
following result.

\begin{proposition}\label{dlbp2}
Let $X$ be a Banach space satisfying $(H_{\beta})$ (resp. $(H^*_{\beta})$)
and $E$ a closed subset of  $X.$ For an upper semicontinuous and bounded from above function $f$ on $E$ and any $\varepsilon \in (0,1),$ there exist a $g\in \mathcal{D}_{\beta}(X)$ (resp. $g\in\mathcal{D}^*_{\beta}(X)$) and an $x_0\in E$ such that:
\begin{enumerate}[(a)]
\item $f-g$ attains its maximum at $x_0.$
\item $\|g\|_{\infty}\le \varepsilon$ and $\|\nabla_{\beta}g\|_{\infty}\le \varepsilon.$
\end{enumerate}
\end{proposition}

\section{The main results}

\begin{definition}
A function $u: \Omega\to\R$ is said to be
\begin{itemize}
\item[(i)] a $\beta$-viscosity subsolution of  \eqref{A1} 
if $u$ is upper semicontinuous and for any $x\in \Omega,$ $x^*\in D^+_{\beta} u(x),$ $ u(x)+H(x, u(x), x^*)\le 0;$
\item[(ii)] a $\beta$-viscosity supersolution of (\ref{A1})
if $u$ is lower semicontinuous and  for any $x\in \Omega,$ $x^*\in D^-_{\beta} u(x),$ $u(x)+H(x, u(x), x^*)\geq 0;$
\item[(iii)] a $\beta$-viscosity solution of  (\ref{A1}) if  $u$ is simultaneously
a $\beta$-viscosity subsolution and a $\beta$-viscosity supersolution.
\end{itemize}
\end{definition}
 
 For convenience, hereafter, we will use the phrases ``$\beta$-viscosity solution of $u+H(x,u,Du)\leq 0$"
and ``$\beta$-viscosity subsolution of $u+H(x,u,Du)=0$" interchangeably. Similarly for the phrases ``$\beta$-viscosity solution of $u+H(x,u,Du)\geq 0$" and ``$\beta$-viscosity supersolution of $u+H(x,u,Du)=0$".

\begin{definition}
A function $u: \Omega\to\R$ is said to be a $\beta$-viscosity subsolution (resp. supersolution, solution) of the problem \eqref{A1}-\eqref{A2} if $u$ is a $\beta$-viscosity subsolution (resp. supersolution, solution) of Equation \eqref{A1} and $u\le \varphi$ (resp. $u\ge \varphi, u=\varphi$) on $\partial \Omega.$
\end{definition}

First, let us recall some auxiliary results. A function $m: [0,\infty) \to  [0,\infty)$
is called a modulus if it is continuous, nonnegative nondecreasing, subadditive on
$[0,\infty)$, and $m(0)=0.$ We will also say that a function
$\sigma: [0,\infty)\times[0,\infty)\to [0,\infty)$ is local modulus if
the function $\sigma(r,R)$ is a modulus for each $R\geq 0$ and $\sigma(r,R)$
is continuous and nondecreasing in both variables.

\begin{remark}\label{nxmd}
It should be noted that if $m: [0,+\infty)\to [0,+\infty)$ is a modulus
then there exist real scalars $A, B$ such that $m(x)\le A+Bx$ for all $x\in [0,+\infty).$
Indeed, for any $x\in [0,+\infty),$ let $[x]$ be the integer part of $x$ and
$A$ is the maximum value of $m$ on $[0,1]$ then $m(x- [x])\leq A.$
Therefore,
\begin{align*}
 m(x)&= m(\underbrace{1+1+\ldots+1}_{ [x] \text{ numbers }}+x- [x])  \\
&\leq [x]m(1)+m(x- [x])\le m(1)x+A=A+Bx,
\end{align*}
where $B=m(1).$
\end{remark}

Next, we make the following assumptions on the function $H.$
\begin{itemize}
\item (H0) There exists a continuous function $w_R: X^*_{\beta} \to \mathbb{R}$ for each $R>0,$ satisfying 
\begin{gather*}
|H(x,r,p)-H(x,r,q)|\le w_R(p-q) 
\end{gather*}
for any $x\in X,$ $p,q\in X^*$ and $r\in \R$ such that $|x|, |q|, |p|\leq R.$
\item (H1) For each $(x,p)\in X\times X^*,$ $r\mapsto H(x,r,p)$ is nondecreasing. 
\item(H1)* For each $(x,p)\in X\times X^*,$ $r\mapsto H(x,r,p)$ is Lipschitz continuous with constant $L_H<1.$ 
\item (H2) There is a local modulus $\sigma_H$ such that
\begin{gather*}
H(x,r,p)-H(x,r,p+q)\le \sigma_H(|q|,|p|+|q|) 
  \end{gather*}
for all $r\in \mathbb{R},$ $x\in\Omega $ and $p, q\in X^*.$
\item (H3) There is a modulus $m_H$ such that
 \begin{equation}
 \begin{split}
H(y,r,\lambda (\nabla_{\beta}|.|^2)(x-y))-&H(x,r,\lambda (\nabla_{\beta}|.|^2)(x-y))\\
&\leq m_H(\lambda |x-y|^2+ |x-y|)
 \end{split}  \notag
\end{equation}
for all $x,y\in \Omega$ with $x\neq y,$ 
$r\in\R$ and $\lambda\geq 0.$
\end{itemize}

\begin{example}
Suppose that $X$ is a Banach space with a $\beta$-smooth norm.
Let $b: X\to X$ be a bounded and Lipschitz continuous function. Consider the following functions: 
$$H(x,r,p)=\langle p,b(x)\rangle \text{ and } H^*(x,r,p)=-\frac{1}{2}r+\langle p,b(x)\rangle.$$
It can be verified that the function $H$ satisfies hypotheses (H0), (H1), (H2) and (H3). Specifically, if $b$ is bounded on $X$ by a constant $K,$ and its Lipschitz constant is $L_b,$ then
\begin{align*}
|H(x,r,p)-H(x,r,q)|&=|\langle p-q,b(x)\rangle| \le |p-q|.|b(x)|\\
&\le K|p-q|=: w_R(p-q)
\end{align*}
for any $x\in X,$ $p,q\in X^*,r\in \R$ satisfying $|x|, |q|, |p|\le R.$
Thus, (H0) is satisfied. Moreover, (H1) and (H2) are obvious. We also have
 \begin{align*}
  H(y,r,\lambda (\nabla_{\beta}|.|^2)(x-y))&-H(x,r,\lambda (\nabla_{\beta}|.|^2)(x-y))\\
 &=\langle\lambda (\nabla_{\beta}|.|^2)(x-y),b(y)-b(x)\rangle\\
&\leq 4KL_b\lambda |x-y|^2\\
&\leq 4KL_b(\lambda |x-y|^2+|x-y|)=: m_H(\lambda |x-y|^2+ |x-y|).
 \end{align*}
Therefore, (H3) is fulfilled.

Similarly, the function $H^*$ satisfies the hypotheses (H0), (H2), (H3). It should be emphasized that the function $H^*$ does not satisfy  (H1) while it satisfies $(H1)^*$. 
\end{example}

On the basis of the above preliminaries, we are now in a position to present our main result on the uniqueness of $\beta$-viscosity solutions
of \eqref{A1}.

\begin{theorem}\label{dl3.1}
Let $X$ be a Banach space with a $\beta$-smooth norm, and $\Omega$ a open subset of $X.$
Assume that the function $H$ satisfies assumptions (H0), (H1) (resp. (H1)*), (H2), (H3), and  
$\widehat{H}$ satisfies (H0). Let $u,v\in C(\overline{\Omega})$ respectively be
$\beta$-viscosity solutions of the problems
\begin{equation}
u+H(x,u,Du)\leq 0  \text{ and } v+\widehat{H}(x,v,Dv)\ge 0 \text{  on }\Omega,
\label{1.9}
\end{equation}
and suppose that there is a modulus $m$ such that
\begin{equation}\label{1.10}
 |u(x)-u(y)|+ |v(x)-v(y)|\le m(|x-y|)\text{ in } \Omega.
\end{equation}
Then we have
\begin{equation}\label{1.11}
\begin{split} 
  u(x)-v(x)&\leq \sup_{\partial \Omega} (u-v)^+
  +\sup_{\Omega\times \R\times X^*}(\widehat{H}-H)^+,\\
 \Big(resp.\quad  u(x)-v(x)&\leq \sup_{\partial \Omega} (u-v)^+
  +\frac{1}{1-L_H}\sup_{\Omega\times \R\times X^*}(\widehat{H}-H)^+\Big),
\end{split}
\end{equation}
for all $x\in \Omega.$ 

In particular, when $\Omega = X,$ we have the estimate \eqref{1.11} in which the term $\sup_{\partial \Omega} (u-v)^+$ on the right hand side is replaced by zero.
\end{theorem}

The proof of this theorem is rather long. To
facilitate the reading, it is postponed to Section 4. 

\begin{corollary}[Comparison and uniqueness]\label{hq3.1}
Given $X$ a Banach space, and equipped with a $\beta$-smooth norm. Let $\Omega\subset X$ be an open set with boundary $\partial\Omega\not=\emptyset,$ $\varphi$ a continuous function on $\partial \Omega.$ Assume that the function $H$ satisfies the assumptions (H0), (H1) (resp. (H1)*), (H2), and (H3). If $u,v\in C(\overline{\Omega})$ respectively are 
$\beta$-viscosity subsolution and $\beta$-viscosity supersolution of Equation \eqref{A1} satisfying \eqref{1.10}, then $u\leq v$ in $\Omega,$ provided that $u\leq v$ on $\partial \Omega.$ Therefore, the problem \eqref{A1}, \eqref{A2} has at most a solution in $C(\overline{\Omega}).$ 

In the case $\Omega$ is the whole space $X,$ the comparison and the uniqueness of the solution for Equation \eqref{A1} is an obvious consequence. 
 \end{corollary}
\begin{proof}
The comparison follows immediately from Theorem \ref{dl3.1}. Let $u,v\in C(\overline{\Omega})$ be the solutions of the problem \eqref{A1}, \eqref{A2}. Then, $u$ is a subsolution, $v$ is a supersolution of \eqref{A1} and $u\leq v$ on $\partial \Omega,$ so by the comparison result, $u\leq v$ in $\Omega.$ By changing the role of $u$ and $v,$ we obtain that $u\geq v,$ and therefore $u=v$ in $\Omega.$ \hfill $\square$  
\end{proof}

\begin{proposition}\label{md3.1}
If in Corollary \ref{hq3.1}, the space $X$ satisfies either $(H_{\beta})$ or $(H^*_{\beta})$ but it has neither $\beta$-smooth norm  nor  norm being equivalent to a $\beta$-smooth norm, then the conclusions of  Corollary \ref{hq3.1} are not true. 
\end{proposition}
\begin{proof}
Indeed, let $H(x,u,Du)=1+u$ in Equation \eqref{A1}. It is easy to  see that the function $u\equiv -\dfrac{1}{2}$ is a $\beta$-viscosity solution of Equation \eqref{A1}. We will show that $u=|x|$ is another solution. 
Indeed,  we take the space $X$ satisfying $(H_{\beta})$ and having a norm which is not $\beta$-differentiable at any point (such a space $X$ can be the one that Remark II.9, \cite{DGZ93} has pointed out). Then $u$ is obviously  a $\beta$-viscosity supersolution of Equation \eqref{A1}. On the other hand, since $u=|x|$  is a convex function and not $\beta$-differentiable at any point on  $X,$  $D^+_{\beta}u(x)=\emptyset$ for all $x\in X.$ Thus $u$ is a  $\beta$-viscosity subsolution of equation \eqref{A1}.
\hfill $\square$
\end{proof}

\begin{remark}\label{can1}
	
	\begin{enumerate}[1)]
		\item When  $\beta =F$ is Fr\'{e}chet borno, we get the result of Crandall M. G. and Lions P. L. in \cite{CL85}.
		
		\item If  the space $X$ has no $\beta$-smooth norm  but it has a  norm which is equivalent to a $\beta$-smooth norm, we also have the conclusion of  Theorem \ref{dl3.1} by using the equivalent norm instead of the original one.
		
		\item In \cite[Theorem 3.6]{leduc}, the authors confirmed that the hypothesis of existing a Fr\'{e}chet-smooth bump function   is sufficient for the existence of an equivalent Fr\'{e}chet-smooth norm on space $X.$ Thus, Theorem \ref{dl3.1}  is a generalization of   Crandall's result in \cite{dau}. Note that in this work, we do not use the  Radon-Nikodym property of the space $X.$
		\item By \cite[page 211]{DGZ93}, the assumptions that there exists a  function $d: X\to \mathbb{R}$ which is bounded and Fréchet differentiable on $X\backslash \{0\},$ and there is a real number $r>0$ such that $d(x)\ge r |x|$ used in \cite{dau} are equivalent to the existence of a Lipschitz and Fr\'{e}chet differentiable bump function on $X.$ Hence, if there is a such function $d,$ then Theorem \ref{dl3.1} still holds without assumption on the existence of a Fr\'{e}chet-smooth norm or the Radon-Nikodym property as in \cite{dau}.	
	\end{enumerate} 
\end{remark}
\begin{proposition}\label{md3.2}
Let $X$ be a normed space  without a norm being equivalent to a $\beta$-smooth norm but satisfies the assumption $(H_{\beta}).$ Then by adding the boundedness of the involving solutions, Theorem \ref{dl3.1} is still valid.
\end{proposition}
\begin{proof}
  By \cite[Lemma 2.15]{BZ96},  there exist a function $d: X\to \mathbb{R^+}$  and  $K>1$ such that
\begin{enumerate}[i)]
	\item $d$ is bounded, Lipschitz continuous on $X$, and $\beta$-smooth on $X\setminus \{0\};$
	\item $|x|\le d(x)\le K |x|$ if $|x|\le 1$ and $d(x)=2$ if $|x|\ge 1.$
\end{enumerate}

Now assume that $u,v$ are bounded on $\Omega.$ In the proof of Theorem \ref{dl3.1}, we replace $|.|^2$ by $d^2(\cdot)$  and consider the function
\begin{equation*}
\Phi(x,y)=u(x)-v(y)-\left(\frac{d(x-y)^2}{\varepsilon}+\lambda\zeta(|x|-R)\right)
\label{1.161}
\end{equation*}
The conclusion follows by using a similar argument as the proof of Theorem \ref{dl3.1}.
\hfill $\square$
\end{proof}

We emphasize that the hypothesis on the boundedness of  $u,v$ in Proposition \ref{md3.2} is necessary. A  counterexample can be found in the proof of Proposition \ref{md3.1}.


 We proceed to study the stability of the $\beta$-viscosity solutions. Using this stability in the same way as in \cite{DGZ93}, we obtain Proposition \ref{dlhb}.

\begin{theorem}[Stability]\label{dl3.2}
	Let $X$ be a Banach space with a $\beta$-smooth norm, and $\Omega$ an  open subset of $X.$
	Let $u_n\in C(\Omega)$ and $H_n\in C(\Omega\times \mathbb{R}\times X^*_{\beta}),$ $n=1, 2, ...$ converge to $u, H$ respectively as $n\to \infty$ in the following way:
	
	For every $x\in \Omega$ there is an $R>0$ such that $u_n\to u$ uniformly on $B_R(x)$ as $n\to\infty,$  and if $(x,r,p), (x_n,r_n,p_n)\in \Omega\times \mathbb{R}\times X^*_{\beta}$ for $n=1, 2, ...$ and $(x_n,r_n,p_n) \to (x,r,p)$ as $n\to\infty,$ then $H_n(x_n,r_n,p_n)\to H(x,r,p).$
	If $u_n$ is a $\beta$-viscosity supersolution (respectively, subsolution) of $H_n=0$ in $\Omega,$ then $u$
	is a $\beta$-viscosity supersolution (respectively, subsolution) of $H=0$ in $\Omega.$
\end{theorem}
\begin{proof}
	We treat the case that the $u_n$ are $\beta$-viscosity supersolutions of $H_n= 0;$ the case of $\beta$-viscosity subsolutions is entirely similar. Let $x\in \Omega$ and $p\in D^-_{\beta} u(x),$ we can find a local Lipschitz continuous and $\beta$-differentiable function $\varphi: X\to \mathbb{R}$  on $D(\varphi)$ such that:
	\begin{itemize}
		\item[(i)] $(u-\varphi)(x)=0$ and $(u-\varphi)(y)\ge 0$ if $y\in \Omega.$
		\item[(ii)] $p=\nabla_{\beta}\varphi(x)$ and $\nabla_{\beta}\varphi$ is norm to $\tau_{\beta}$ continuous at $x.$
	\end{itemize}
	Since $u_n$ converges to $u$ as $n\to \infty$, then we can assume that there exist a sequence of positive numbers $\{\varepsilon_n\}_n$  such that $\varepsilon_n \to 0$ 
	as $n\to \infty$ and
	$$u(y)\le u_n(y)+\varepsilon_n \text{ for }y\in \Omega.$$
	Therefore,  $$u_n(y)+\varepsilon_n-\varphi(y)\ge(u-\varphi)(y)\ge 0\text{ for }y\in \Omega,$$
	which implies $\inf (u_n-\varphi)+\varepsilon_n\ge 0.$
	Since $u$ and $\varphi$ are continuous at $x,$ there exist $z_n\in \Omega,$ $|z_n-x|<\varepsilon_n$ such that 
	$u_n(z_n)-\varphi(z_n)\le \varepsilon_n$. Thus,  $u_n(z_n)-\varphi(z_n)\le \inf (u_n-\varphi)+2 \varepsilon_n.$ Let us define $(u_n-\varphi)(y)=-\varepsilon_n$ for $y\notin \Omega$. Then $u_n-\varphi$ is lower semicontinuous and bounded from  below. According to Proposition \ref{dlbp}, there exists $\psi_n\in \mathcal{D}_{\beta}(X)$ such that 
	\begin{itemize}
		\item[(1)] $\lim\limits_{n\to \infty}||\psi_n||=0$ and $\lim\limits_{n\to \infty}||\nabla_{\beta}\psi_n||_{\infty}=0,$
		\item[(2)] $u_n-\varphi+\psi_n$ attains a local minimum at $y_n\in \Omega,$
		\item[(3)] $\lim\limits_{n\to \infty}|y_n-z_n|=0.$
	\end{itemize}
	Using (3) and the fact that $|z_n-x|<\varepsilon_n,$ we have $\lim\limits_{n\to \infty}|y_n-x|=0.$ On the other hand, $p_n=\nabla_{\beta}\varphi(y_n)-\nabla_{\beta}\psi_n(y_n)\in D^-_{\beta} u_n(y_n) $ and $p_n-p=\nabla_{\beta}\varphi(y_n)-\nabla_{\beta}\varphi(x) -\nabla_{\beta}\psi_n(y_n).$
	Since $\nabla_{\beta}\varphi$ is a continuous function, $y_n\to x$ and $\lim\limits_{n\to \infty}||\nabla_{\beta}\psi_n||_{\infty}=0,$ $p_n-p\to 0$ when $n\to \infty.$
	
	Now using the fact that $u_n$ is a $\beta$-viscosity supersolution of $H_n=0,$ we obtain 
	$$H_n(y_n,u_n(y_n),p_n)\ge 0.$$
	Passing to the limit as $n\to \infty$, we get $H(x,u(x),p)\ge 0.$ Thus, $u$ is $\beta$-viscosity supersolution of $H(x,r,p)=0.$
	\hfill $\square$  
\end{proof}

Finally, we close this section with a finding on the existence of $\beta$-viscosity solutions. The following Theorem \ref{dl3.3} generalizes the ones in \cite{CL86} and \cite{DGZ93}. To prepare for that theorem, we recall several definitions and a proposition proved in \cite{DGZ93}.

If $\Omega$ is an open subset of a Banach space $X$ and if $u$ is a function defined on $\Omega,$   the upper semicontinuous envelope $u^*$ of $u$ is defined as
$$u^*=\inf\{v: v\text{ is continuous on }\Omega\text{ and }v\ge u\text{ on }\Omega\}$$
and the lower semicontinuous envelope $u_*$ of $u$ is defined by
$$u_*=\sup\{v: v\text{ is continuous on }\Omega\text{ and }v\le u\text{ on }\Omega\}.$$

\begin{proposition}[Proposition III.7, \cite{DGZ93}]\label{dlhb}
	Let $X$ be a Banach space satisfying $(H^*_{\beta}),$  $\Omega$ an open subset of $X,$  and $F: \overline{\Omega} \times\mathbb{R} \times X^*_\beta \to \mathbb{R}$ continuous. Suppose that $u_0, v_0$ are  respectively a $\beta$-viscosity subsolution and a $\beta$-viscosity supersolution of $F(x,u,Du)=0$ on $\Omega.$ Assume $u_0\le v_0$ on $\Omega.$
	Then, there is a function $u: X\to \mathbb{R}$ such that $u_0\le u\le v_0$ on $\Omega,$ $u^*$ is a $\beta$-viscosity subsolution of $F(x,u,Du)=0$ on $\Omega,$ and $u_*$ is a $\beta$-viscosity supersolution of $F(x,u,Du)=0$ on $\Omega.$
\end{proposition}

We are now ready to establish the existence.

\begin{theorem}[Existence]\label{dl3.3}
	Let $X$ be a Banach space with a $\beta$-smooth norm, and $\Omega$ an  open subset of $X.$ Let $H: \Omega \times\mathbb{R} \times X^* \to \mathbb{R}$ satisfy (H0), (H1) (resp. (H1)*), (H2), (H3), $\varphi\in C(\partial \Omega)$ such that $|\varphi (x)|\le M+N|x|$ for some real numbers $M, N\in (0,\infty)$ and for all $x\in \partial\Omega.$ Then there exists a unique $\beta$-viscosity solution of the problem \eqref{A1}, \eqref{A2}.
\end{theorem}
\begin{proof}
	The uniqueness part of the problem is asserted in	Corollary \ref{hq3.1}. Here we only present the proof of the existence using the Perron's method. 
	
	By taking $\lambda=0$ and $r=0,$ the assumption (H3) gives us $H(x,0,0)- H(0,0,0)\le m_H(|x|),$ so
	$$ H(x,0,0)\le H(0,0,0)+m_H(|x|) \le |H(0,0,0)|+m_H(|x|). $$
Similarly, we have $- H(x,0,0)\le |H(0,0,0)|+m_H(|x|).$ Therefore $|H(x,0,0)|\le |H(0,0,0)|+m_H(|x|).$

	Since $m_H$ is a modulus, according to Remark \ref{nxmd}, there exist real numbers $A_H, B_H$ such that 
	$$	|H(x,0,0)|\le A_H+B_H |x|,\quad \forall x\in \Omega.$$
	
	Let $B_1= \max\{B_H,N\}.$  Defining  $$A_1=\max\{A_H+\sigma_H(B_H,B_H), M\},\quad  A_2=\max\{A_H+\sigma_H(B_1,2B_1), M\},$$
	$$A'_1=\max\Big\{A_H+\sigma_H\Big(\frac{B_1}{1-L_H},\frac{B_1}{1-L_H}\Big), M\Big\};$$ $$A'_2=\max\Big\{A_H+\sigma_H\Big(\frac{B_1}{1-L_H},\frac{2B_1}{1-L_H}\Big), M\Big\}$$
	and
	$$v_0(x)=A_1+B_1|x|,\quad u_0(x)=-(A_2+B_1|x|).$$
	$$\Big({\rm resp.}\quad v_0(x)=\frac{A'_1+B_1|x|}{1-L_H},\quad u_0(x)=-\frac{A'_2+B_1|x|}{1-L_H}\Big).$$
	
	We claim that $v_0$ is a $\beta$-viscosity supersolution and $u_0$ is a $\beta$-viscosity subsolution of the problem \eqref{A1}, \eqref{A2} in $\Omega.$ 
	Indeed, it is clear that $v_0\geq \varphi$ on $\partial \Omega.$ Moreover, from the assumptions (H1) (resp. (H1)*) and (H2), for each $x\in \Omega,$ $p\in D^-_{\beta} (|\cdot|)(x),$ we have 
	\begin{align*}
	H(x,0,0)&\le H(x,A_1+B_1|x|,0)\\
	& \le H(x,A_1+B_1|x|,B_1p)+\sigma_H(B_1,B_1);
	\end{align*}	
	\begin{align*}\Big({\rm resp.} \quad 
	H(x,0,0)&\le H(x,v_0(x),0)+L_Hv_0(x)\\
	& \le H\Big(x,v_0(x),\frac{B_1}{1-L_H}p\Big)+\sigma_H\Big(\frac{B_1}{1-L_H},\frac{B_1}{1-L_H}\Big)+L_Hv_0(x)
	 \Big). \end{align*}
	Hence
	\[\begin{split}
	v_0(x)+H(x, v_0(x),B_1p)&\ge A_1+B_1|x|+ H(x,0,0)-\sigma_H(B_1,B_1)\\
	&\ge A_1+B_1|x|- A_H-B_H|x|-\sigma_H(B_1,B_1)\\
	&\ge 0;
	\end{split}\]
	\[\begin{split}
	\Big({\rm resp.}\quad  v_0(x)+&H\Big(x, v_0(x),\frac{B_1}{1-L_H}p\Big)\\
	&\ge (1-L_H)v_0(x)+ H(x,0,0)-\sigma_H\Big(\frac{B_1}{1-L_H},\frac{B_1}{1-L_H}\Big)\\
	&\ge A'_1+B_1|x|- A_H-B_H|x|-\sigma_H\Big(\frac{B_1}{1-L_H},\frac{B_1}{1-L_H}\Big)\ge 0\Big).
	\end{split}\] 
	This proves that $v_0$ is a $\beta$-viscosity supersolution of \eqref{A1}. 
	
	On the other hand, by the definition, $u_0\leq\varphi$ on $\partial\Omega.$ It remains to prove that $u_0$ is a $\beta$-viscosity subsolution of \eqref{A1} in $\Omega.$
	
	From the conditions (H1) (resp. (H1)*) and (H2), for each $x\in \Omega,$ $p\in D^+_{\beta}(|~\cdot~|)(x),$ we have 
	\begin{gather*}
	H(x,u_0(x),-B_1p)\le H(x,0,-B_1p)\le H(x,0,0)+\sigma_H(B_1,2B_1).
	\end{gather*}
	\[\begin{split}\Big({\rm resp.}\quad 
	H\Big(x,u_0(x),&-\frac{B_1}{1-L_H}p\Big)\le H\Big(x,0,-\frac{B_1}{1-L_H}p\Big)+L_H|u_0(x)|\\
	&\le H(x,0,0)+\sigma_H\Big(\frac{B_1}{1-L_H},\frac{2B_1}{1-L_H}\Big)+L_H|u_0(x)|\Big).
	\end{split}\]
	Hence
	\[\begin{split}
	u_0(x)+H(x,u_0(x),-B_1p)&\le u_0(x)+H(x,0,0)+\sigma_H(B_1,2B_1)\\
	&\le u_0(x)+A_H+B_H|x|+\sigma_H(B_1,2B_1)\le 0;
	\end{split}\]
	\[\begin{split} 	\Big({\rm resp.}\quad 
	u_0(x)+&H\Big(x,u_0(x),-\frac{B_1}{1-L_H}p\Big)\\
	&\le (1-L_H)u_0(x)+H(x,0,0)+\sigma_H\Big(\frac{B_1}{1-L_H},\frac{2B_1}{1-L_H}\Big)\\
	&= -(A'_2+B_1|x|)+A_H+B_H|x|+\sigma_H\Big(\frac{B_1}{1-L_H},\frac{2B_1}{1-L_H}\Big)\le 0\Big).
	\end{split}\]
	Thus $u_0$ is a $\beta$-viscosity subsolution of \eqref{A1}.
	
	Now by Proposition \ref{dlhb}, there exists a function $u$ such that $u_0\le u \le v_0,$ $ u_*$ and $u^*$ are respectively a $\beta$-viscosity supersolution and a $\beta$-viscosity subsolution of Equation \eqref{A1} in $\Omega.$ It follows from Corollary \ref{hq3.1} that $u_*\ge u^*.$ Clearly by the definition,  $u_*\le u^*.$ Hence $u^*=u_*=u$ is a $\beta$-viscosity solution of Equation \eqref{A1} in $\Omega.$ 	\hfill $\square$  	
\end{proof}
\begin{example}\label{vd0T1}
	In Example 3.1, we look at
	a linear Hamiltonian $H$ and verify that assumptions (H0)-(H3) are satisfied. Here is another example, in which the Hamilton function $H$ is nonlinear. Moreover, the underlying normed space $X$ below has no equivalent Fr\'{e}chet smooth norm.

	Let $X=L_1[0,1].$ By \cite{BZ96}, $X$ has an equivalent weak Hadamard smooth norm and that norm is nowhere Fr\'{e}chet-differentiable. 
	Thus, the results in Crandall \cite{CL85} do not work in this case. Also note that $L_1[0,1]$ does not satisfy the Radon-Nikodym property  (see \cite{BZ99}). Let $U=[0,1]$ and $\mathcal{U}=\{ \alpha: [0,\infty) \to U \text{ is continuous}\}$. Let a function $H: X\times X_{\beta}^*\to \mathbb{R}$ defined by   
	$$H(x,p)=\sup_{\alpha \in \mathcal{U}} \{-\langle p, x_0 \alpha \rangle -|x| \alpha  \},$$
	where $x_0$ is a fixed point in $X$. Then it can be verified that $H$ satisfies conditions (H0)-(H3). 
	Consider the equation 
	\begin{equation}\label{TPVD}
	u+H(x,Du)=0 \quad \text{ on } X.
	\end{equation}
	By virtue of Proposition \ref{md3.2} and Theorem \ref{dl3.2}, equation \eqref{TPVD} has a unique weak Hadamard-viscosity solution.
	\end{example}

\section{Proof of Theorem \ref{dl3.1}}
\subsection{Proof for the case $\Omega\not= X$}
 The conclusion in this case is clear if $\sup_{\partial  \Omega}(u-v)^+=\infty.$ Thus, we now assume  $\sup_{\partial  \Omega}(u-v)^+<\infty.$

To continue, we first give the following lemma. 
\begin{lemma}\label{bd3.1}
	Suppose that $\Omega\subset X$ is an open strict subset. We denote as $\rho(x)$ the distance from $x$ to $\partial\Omega.$ Let $u,v\in C(\overline{\Omega})$ be functions that satisfying the condition (\ref{1.10}). Then, for all $x,y \in  \Omega,$ 
	\begin{equation}\label{1.13}
	u(x)-v(y)\le \sup_{\partial  \Omega}(u-v)^+ + m(\min(\rho(x),\rho(y)))+m(|x-y|).
	\end{equation}
	In particularly, 
	\begin{equation}\label{T3.7}
	u(x)-v(x) \le \sup_{\partial  \Omega}(u-v)^+ + m(\rho(x)), \ x\in \Omega.
	\end{equation} 
\end{lemma}
\begin{proof}
	Suppose that $\rho(x)< \rho(y).$ Then, for any $\varepsilon>0,$ there exists
	an $x^*\in \partial \Omega$ such that  the line segment connecting $x$ and $x^*$ is contained in $\Omega,$ and $|x-x^*|<\rho(x)+\varepsilon.$ 
	For each $n\in \mathbb{N^*},$ we choose $x_n=x+\left(1-\frac{1}{n}\right)(x^*-x).$ Then, 
	$x_n\in\Omega,$ and $x_n\to x^*$ as $n\to \infty.$ In addition,
	\begin{align}
	u(x)-v(y)&=[u(x)-u(x_n)]+[v(x_n)-v(x)] +[u(x_n)-v(x_n)] +[v(x)-v(y)] \notag \\
	&\le m(|x-x_n|)+u(x_n)-v(x_n)+m(|x-y|). \label{them3}
	\end{align}
	Let $n\to \infty$ in \eqref{them3}, we obtain
	\begin{equation*}
	u(x)-v(y)\le m(|x-x^*|)+u(x^*)-v(x^*)+m(|x-y|).
	\label{them4}
	\end{equation*}
	Moreover, $u(x^*)-v(x^*)\le \sup_{\partial  \Omega}(u-v)^+ $ and
	$m(|x-x^*|)\le m(\rho(x)+\varepsilon).$ Thus,
	\begin{equation}
	u(x)-v(y)\le m(\rho(x)+\varepsilon)+\sup_{\partial  \Omega}(u-v)^+ +m(|x-y|).
	\label{them5}
	\end{equation}
	The inequality \eqref{them5} holds for any $\varepsilon>0.$
	By letting $\varepsilon\to 0$ and using the continuity of the function $m$ at $\rho(x),$ we obtain
	$$ u(x)-v(y)\leq \sup_{\partial  \Omega}(u-v)^+ + m(\rho(x))+m(|x-y|). $$
	If  $\rho(x)\geq \rho(y),$ by a similar argument, we have
	\begin{equation*}
	u(x)-v(y)\le \sup_{\partial  \Omega}(u-v)^+ + m(\rho(y))+m(|x-y|).
	\label{them8}
	\end{equation*}
	Thus, \eqref{1.13} is proved.\hfill $\square$
\end{proof}

By Lemma \ref{bd3.1}, the estimate \eqref{T3.7}, and the function 
$\rho(x)$ has at most linear growth, there are constants $A,B \geq 0$ such that
\begin{equation}\label{1.14}
u(x)-v(x)\le A+B|x|, \ x\in \Omega.
\end{equation}

Next, let $\zeta\in C^1(\R)$  satisfy
	\begin{equation*}
	\zeta(r)=0  \text{ with } r\le 1, \quad \zeta(r)=r-2 \text{ with }r\ge 3,\quad  \text{ and } 0\le \zeta'\le 1.
	\label{1.15}
	\end{equation*}
	For $a, \varepsilon, \lambda, R>0,$ consider the function
	\begin{equation}
	\Phi(x,y)=u(x)-v(y)-\left(\frac{|x-y|^2}{\varepsilon}+\lambda\zeta(|x|-R)\right)
	\label{1.16}
	\end{equation}
	on the set
	\begin{equation*}
	\Delta(a)=\{(x,y)\in \Omega\times \Omega: \rho(x),\rho(y)>a \text{\ \  and\ \  } |x-y|<a\}.
	\label{1.17}
	\end{equation*}
We claim that for $B$ defined in \eqref{1.14}, if 
	\begin{equation*}
	\lambda >B,\quad \varepsilon\le a^2/(m(a)+1)\quad \text{  and  }\quad R>1,
	\label{1.18}
	\end{equation*}
	then
	\begin{equation}\label{1.19}
	\begin{split}
	\Phi(x,y)&\le \sup_{\partial \Omega} (u-v)^++2m(a)+C(a,\lambda,\varepsilon)\\
	\Big({\rm resp.}\quad \Phi(x,y)&\le \sup_{\partial \Omega} (u-v)^++2m(a)+ \frac{1}{1-L_H}C(a,\lambda,\varepsilon)\Big)
	\end{split}
	\end{equation}
	on $\Delta(a),$ where
	\begin{equation*}
	\begin{split}
	C(a,\lambda,&\varepsilon)=m_H(2m(a)+(\varepsilon m(a))^{1/2})+\sigma_H(\lambda,2(m(a)/\varepsilon)^{1/2}+\lambda)\\
	&+\sup\{(\widehat{H}-H)^+(z,r,p): (z,r,p)\in\Omega\times \R\times X^*\text{ and }|p|\le 2(m(a)/{\varepsilon})^{1/2}\}
	\end{split}
	\end{equation*}
	The proof of \eqref{1.19} is carried out in $2$ steps.
	
	\medskip 
	\noindent 
	\textbf{Steps 1.} Proof of \eqref{1.19} in the case $ \Phi$ attains the maximum value on $\Delta(a)$ at point $(x_0,y_0).$
	
	If it is the case, fix $y=y_0.$ Then the function $u(x)-\left(\frac{|x-y_0|^2}{\varepsilon}+\lambda\zeta(|x|-R)\right)$ attains a local maximum at $x=x_0.$ Hence 
	$$\frac{1}{\varepsilon}(\nabla_{\beta}|.|^2) (x_0-y_0)+\lambda\zeta'(|x_0|-R).(\nabla_{\beta}|.|)(x_0)    \in D^+_{\beta}u(x_0).$$
	Similarly, when we fix $x=x_0,$ the function $v$ is $\beta$-viscosity subdifferentiable at $y_0$ and
	$$\frac{1}{\varepsilon}(\nabla_{\beta}|.|^2) (x_0-y_0)    \in D^-_{\beta}v(y_0).$$
	Since $u,v$ are $\beta$-viscosity solutions of equation  (\ref{1.9}), we have
	\begin{align}
	u(x_0)+H(x_0,u(x_0), p_{\varepsilon}+\lambda q)&\le 0,\label{1.311} \\
	v(y_0)+\widehat{H}(y_0,v(y_0),p_{\varepsilon})&\ge 0,
	\label{1.312}
	\end{align}
	where $p_{\varepsilon} =\frac{1}{\varepsilon}(\nabla_{\beta}|.|^2) (x_0-y_0), q=\zeta'(|x_0|-R).(\nabla_{\beta}|.|)(x_0).$

	Because the function  $\Phi$  attains its maximum value  at $(x_0,y_0)$  so for all $(x,y)\in \Delta(a),$ then
	\begin{equation}
	\Phi(x,y) \le \Phi(x_0,y_0)\le u(x_0)-v(y_0).
	\label{1.191}
	\end{equation}
	
	If $u(x_0)-v(y_0)\le 0,$ then $\Phi(x,y)\leq 0$ in $\Delta(a).$ This fact proves \eqref{1.19} because its right-hand side is nonnegative.
	
	If $u(x_0)-v(y_0)> 0,$ then $u(x_0)>v(y_0).$ From (\ref{1.311}) and (\ref{1.312}), we imply that 
	\begin{equation}
	u(x_0)-v(y_0)\le  \widehat{H}(y_0,v(y_0),p_{\varepsilon})-H(x_0,u(x_0),p_{\varepsilon}+\lambda q).
	\label{1.33} 
	\end{equation}
In addition, using (H1) (resp. (H1)*) we obtain
\begin{equation}\label{1.35}
\begin{split}
\widehat{H}&(y_0,v(y_0),p_{\varepsilon})-H(x_0,u(x_0),p_{\varepsilon}+\lambda q)\\
&= \widehat{H}(y_0,v(y_0),p_{\varepsilon})-H(y_0,v(y_0),p_{\varepsilon})+H(y_0,v(y_0),p_{\varepsilon})-H(x_0,u(x_0),p_{\varepsilon}+\lambda q) \\
&\le  \widehat{H}(y_0,v(y_0),p_{\varepsilon})-H(y_0,v(y_0),p_{\varepsilon})+H(y_0,u(x_0),p_{\varepsilon})-H(x_0,u(x_0),p_{\varepsilon}+\lambda q)\\
&= \widehat{H}(y_0,v(y_0),p_{\varepsilon})-H(y_0,v(y_0),p_{\varepsilon})
+H(y_0,u(x_0),p_{\varepsilon})-H(x_0,u(x_0),p_{\varepsilon})\\
&\ \  +H(x_0,u(x_0),p_{\varepsilon})-H(x_0,u(x_0),p_{\varepsilon}+\lambda q);\\
\Big(&{\rm resp.}\quad 
\widehat{H}(y_0,v(y_0),p_{\varepsilon})-H(x_0,u(x_0),p_{\varepsilon}+\lambda q)\\
&\leq  \widehat{H}(y_0,v(y_0),p_{\varepsilon})-H(y_0,v(y_0),p_{\varepsilon})
+H(y_0,u(x_0),p_{\varepsilon})-H(x_0,u(x_0),p_{\varepsilon})\\
&\ \  +H(x_0,u(x_0),p_{\varepsilon})-H(x_0,u(x_0),p_{\varepsilon}+\lambda q)+L_H(u(x_0)-v(y_0)\Big).
\end{split}
\end{equation}	
	
On the other hand, we have
	\begin{align*}
	0\le \Phi(x_0,y_0)- \Phi(x_0,x_0)&=v(x_0)-v(y_0)-\frac{|x_0-y_0|^2}{\varepsilon} \\ 
	&\le m(|x_0-y_0|)-\frac{|x_0-y_0|^2}{\varepsilon}
	\end{align*}
	so 
	$$ \frac{|x_0-y_0|^2}{\varepsilon} \le  m(|x_0-y_0|) \le m(a),$$
	or $$|x_0-y_0|\le (\varepsilon m(a))^{1/2}.  $$
	
	Together with the assumption (H3) and the increasing monotonicity of the module $m_H,$ we have
	\begin{equation}\label{4.101}
	\begin{split} 
	H(y_0,u(x_0),p_{\varepsilon})-H(x_0,u(x_0),p_{\varepsilon}) &\le m_H(2|x_0-y_0|^2/\varepsilon+|x_0-y_0|)\\ 
	&\le m_H(2 m(a)+(\varepsilon m(a))^{1/2}).
	\end{split} 
	\end{equation}
Because of the facts that
	$$  \lambda |q|=\lambda |\zeta'(|x_0|-R).(\nabla_{\beta}|.|)(x_0)|\le \lambda,$$
	and 
	$$ |p_{\varepsilon}| =\frac{1}{\varepsilon}|(\nabla_{\beta}|.|^2) (x_0-y_0) |\le \frac{2}{\varepsilon}  |x_0-y_0|\le 2 (m(a)/\varepsilon)^{1/2},$$
by using (H2) we yield	
	\begin{equation}\label{4.102}
\begin{split} 
H(x_0,u(x_0),p_{\varepsilon})-H(x_0,u(x_0),p_{\varepsilon}+\lambda q)&\le \sigma_H(\lambda |q|,|p_{\varepsilon}|+\lambda |q|)\\
	&\le \sigma_H(\lambda,2 (m(a)/\varepsilon)^{1/2}+\lambda).
	\end{split} 
\end{equation}

Plugging \eqref{4.101}, \eqref{4.102} into \eqref{1.35}, then combining with \eqref{1.191}, \eqref{1.33} we have
\[\Phi(x,y)\leq C(a,\lambda,\varepsilon)\quad (resp.\quad  \Phi(x,y)\leq \frac{1}{1-L_H}C(a,\lambda,\varepsilon)), \] 
in $\Delta(a)$ and \eqref{1.19} holds in this case.
	
	\medskip 
	\noindent 
	\textbf{Steps 2.} 	Prove (\ref{1.19}) in the case the function $\Phi$ does not attain its maximum value on $\Delta(a).$ 
	
	To do this, firstly we claim that function $\Phi$ is bounded from above on $\Delta(a).$ Indeed, since $\lambda>B,$ for any $x, y\in\Delta(a)$ we have
	\begin{align}
	\Phi(x,y)& \le u(x)-v(x)+v(x)-v(y)-\frac{|x-y|^2}{\varepsilon}-\lambda(|x|-R-2) \notag \\
	&\le A+B|x|+m(a)-\lambda|x|+\lambda(R+2)\to -\infty \text{ khi  }|x|\to \infty. \label{them2}
	\end{align} 
	Then  $\Phi<0$ when $|x|>R_0$ for $R_0$ large enough. On the other hand, for $|x|\le R_0$ we have
	\begin{align*}
	\Phi(x,y)&\le u(x)-v(y)=u(x)-v(x)+v(x)-v(y)\le u(x)-v(x)+m(a)\\
	& \le A+B|x|+m(a)\le A+B R_0+m(a)<\infty
	\end{align*} 
	so  $\Phi(x,y)$ is bounded from above on $\Delta(a).$

	Since $\Phi$ is bounded from above, there exists the finite $\sup_{\Delta(a)}\Phi.$ 
	
	If  $\sup_{\Delta(a)}\Phi\le\sup_{\partial \Omega} (u-v)^++2m(a),$ then (\ref{1.19}) is obvious. Otherwise, we have 
	\begin{equation}
	\sup_{\Delta(a)}\Phi>\sup_{\partial \Omega} (u-v)^++2m(a).
	\label{1.20}
	\end{equation}
		From  (\ref{1.13}) and  (\ref{1.16})   we have
	\begin{equation}
	\Phi(x,y)\le  \sup_{\partial  \Omega}(u-v)^+ + m(\min(\rho(x),\rho(y)))+m(a).
	\label{1.22'}
	\end{equation}
	
	Choose a sequence $(x_n,y_n)\in \Delta(a), n=1,2,...$ such that $\Phi(x_n,y_n)$ is stricly increasing,   
	\begin{equation}
	\Phi(x_n,y_n)\to \sup_{\Delta(a)}\Phi \text{  and  } \Phi(x_n,y_n)\ge \Phi(x_n,x_n).
	\label{1.21}
	\end{equation}
	We prove that $(x_n,y_n)$ is an interior point of $\Delta(a),$  or in other words, there exists $\gamma >0$ such that
	\begin{equation}
	S_n:=\{(x,y)\in X\times X: |x-x_n|^2+|y-y_n|^2\le \gamma^2\}\subset \Delta(a).
	\label{1.26}
	\end{equation}
	Indeed, we have
	\begin{equation}
	\Phi(x,y)-\Phi(x,x)=v(x)-v(y)-\frac{|x-y|^2}{\varepsilon}\le m(|x-y|)-\frac{|x-y|^2}{\varepsilon}
	\label{1.23'}
	\end{equation}
	on $\Delta(a).$ Substituting  $x=x_n,y=y_n$ into (\ref{1.23'}) and combining with (\ref{1.21}), we have
	\begin{equation}
	|x_n-y_n|^2\le \varepsilon m(|x_n-y_n|)\le \varepsilon m(a)\le a^2.m(a)/(m(a)+1).
	\label{1.24}
	\end{equation}
	Hence
	\begin{equation*}
	|x_n-y_n|\le a (m(a)/(m(a)+1))^{1/2}.
	\end{equation*}
	From  (\ref{1.20}) and  (\ref{1.21}) we see that there exists $\xi>0$ such that both 
	\begin{equation}
	\rho(x_n) \text{ and }\rho(y_n)\text{ are not less than } a+\xi \text{ for } n \text{ sufficiently large. }
	\label{1.23}
	\end{equation}
	Indeed, otherwise, we have $\rho(x_n)\to a.$ If $x=x_n,y=y_n$ in (\ref{1.22'})  we have 
	$$  \Phi(x_n,y_n)\le  \sup_{\partial  \Omega}(u-v)^+ + m(\min(\rho(x_n),\rho(y_n)))+m(a).$$
	Then let $n\to \infty$ we obtain $\sup_{\Delta(a)}\Phi\le \sup_{\partial \Omega} (u-v)^++2m(a).$
	This contradicts to (\ref{1.20}). Thus, \eqref{1.23} is proved.
	
	For $(x,y)\in S_n,$ we have $|x-x_n|^2+|y-y_n|^2\le \gamma^2$ so
	$$|x-x_n|+|y-y_n|\le \sqrt{2}. \left(|x-x_n|^2+|y-y_n|^2\right)^{1/2}<2\gamma$$ 
	and $$|x-y|\le |x-x_n|+|x_n-y_n|+|y_n-y|<2\gamma+a (m(a)/(m(a)+1))^{1/2}<a,$$
	where $\gamma$  is a positive number, which is chosen so that 
	$$\gamma<\min\left(\xi, \frac{1}{2}(1-(m(a)/(m(a)+1)))^{1/2}a\right).$$
	Following \eqref{1.23} we have  $\rho(x)\ge \rho(x_n)-|x-x_n|\ge a+\xi -\gamma>a.$ A similar argument can be used to show that $\rho(y)> a,$ hence \eqref{1.26} holds.
	By \eqref{them2} the sequence  $ (x_n,y_n)$  is bounded.

	Next, we will prove that 
	\begin{align}
	\sup_{\Delta(a)}\Phi  &\le C(a,\lambda,\varepsilon)\quad (resp.\quad  \sup_{\Delta(a)}\Phi\leq \frac{1}{1-L_H}C(a,\lambda,\varepsilon)), \label{them1}
	\end{align}
	 which will allow us to obtain  \eqref{1.19}.
	 
	 Indeed, put
	\begin{equation*}
	\delta_n=\sup_{\Delta(a)}\Phi-\Phi(x_n,y_n).
	\label{1.27}
	\end{equation*}
	According to \eqref{1.21}, $\delta_n>0$ for all $n.$
	Consider the function
	\begin{equation}
	\Psi(x,y)=\Phi(x,y)-(3\delta_n/\gamma^2)(|x-x_n|^2+|y-y_n|^2) \text{ on } S_n.
	\label{1.28}
	\end{equation}
	For $(x,y)\in \partial S_n$ we have
	$$ \Psi(x,y) =\Phi(x,y)-3\delta_n\le  \sup_{\Delta(a)}\Phi-3\delta_n=\Phi(x_n,y_n)-2\delta_n=\Psi(x_n,y_n)-2\delta_n.$$
	Hence
	\begin{equation*}
	\Psi(x,y)\le \Psi(x_n,y_n)-2\delta_n\text{ on } \partial S_n.
	\label{1.29}
	\end{equation*}
	Indeed, if $P: S_n\to \R$  varies by less than $\delta_n$ on $S_n$ and $\Psi-P$ has the maximum point in $S_n,$ then this point must belong to the interior of $S_n.$

	By Proposition \ref{dlbp2}, with $\varepsilon =\delta_n, $ there  exists $P\in \mathcal{D}_{\beta}^*(X\times X)$ such that $\Psi(x,y)-P(x,y)$ attains its maximum value on $S_n$ at a point $(\widehat{x},\widehat{y})$ and  $\|P\|_{\infty}<\delta_n,$ $\|\nabla_{\beta}P\|_{\infty}< \delta_n.$ 
	By the argument above, it lies in the interior of  $S_n.$ 
	
	We fix  $y=\widehat{y}.$ Then we have 
		\begin{multline*}
	\Psi(x,\widehat{y})-P(x,\widehat{y})= u(x)-v(\widehat{y})-\left(\frac{|x-\widehat{y}|^2}{\varepsilon}+\lambda\zeta(|x|-R)\right)\\
	-(3\delta_n/\gamma^2)(|x-x_n|^2+|\widehat{y}-y_n|^2) -P(x,\widehat{y})
	\end{multline*}
	attains its maximum on $S_n$ at the point $x=\widehat{x}.$ Hence the function $u$  is $\beta$-viscosity superdifferentiable at $\widehat{x}$ and
	\[\begin{split}
	\frac{1}{\varepsilon}(\nabla_{\beta}|.|^2) (\widehat{x}-\widehat{y})&+\lambda\zeta'(|\widehat{x}-R|)(\nabla_{\beta}|.|)(\widehat{x})\\
	&+(3\delta_n/\gamma^2)(\nabla_{\beta}|.|^2) (\widehat{x}-x_n)+\nabla_{\beta}^xP(\widehat{x},\widehat{y})    \in D^+_{\beta}u(\widehat{x}).
	\end{split}\]
	Similarly,  the function $v$ is $\beta$-viscosity subdifferentiable at $\widehat{y}$ and
	$$\frac{1}{\varepsilon}(\nabla_{\beta}|.|^2) (\widehat{x}-\widehat{y})-(3\delta_n/\gamma^2)(\nabla_{\beta}|.|^2) (\widehat{y}-y_n)-\nabla_{\beta}^yP(\widehat{x},\widehat{y})    \in D^-_{\beta}v(\widehat{y}),$$
	where $\nabla_{\beta}^x$ is the $\beta$-derivative with respect to $x$ variable. 
	
	Since $u,v$ are respectively $\beta$-viscosity subsolution and $\beta$-viscosity supersolution of equation  (\ref{1.9}), 
	\begin{align}
	u(\widehat{x})+H(\widehat{x},u(\widehat{x}), p_{1\varepsilon}+\lambda q+\theta_{1n})&\le 0,\label{1.31'} \\
	v(\widehat{y})+\widehat{H}(\widehat{y},v(\widehat{y}),p_{1\varepsilon}+\theta_{2n})&\ge 0,
	\label{1.31}
	\end{align}
	where
	\begin{equation}
	p_{1\varepsilon}=\frac{1}{\varepsilon}(\nabla_{\beta}|.|^2) (\widehat{x}-\widehat{y}),
	q=\zeta'(|\widehat{x}-R|).(\nabla_{\beta}|.|)(\widehat{x}),
	\label{1.32}
	\end{equation}
	\begin{equation*}
	\theta_{1n}=K_n.(\nabla_{\beta}|.|^2) (\widehat{x}-x_n)+\nabla_{\beta}^xP(\widehat{x},\widehat{y}),\quad \theta_{2n}=-K_n.(\nabla_{\beta}|.|^2) (\widehat{y}-y_n)-\nabla_{\beta}^yP(\widehat{x},\widehat{y}),
	\end{equation*}
	\begin{equation*}
	K_n=3\delta_n/\gamma^2.
	\end{equation*}
	
	Since
	$(\widehat{x},\widehat{y})\in S_n, $  the sequence $(x_n,y_n)$ is bounded, hence $(\widehat{x},\widehat{y})$ and $(x_n,y_n)$  are contained in a bounded set.
	
	Moreover, since  $(\nabla_{\beta}|.|^2)(x)$ is bounded by $2|x|,$ we have $$|p_{1\varepsilon}+\lambda q+\theta_{1n}|\le \frac{2}{\varepsilon}|\widehat{x}-\widehat{y}|+\lambda |\widehat{x}|+K_n |\widehat{x}-x_n|+|\nabla_{\beta}^xP(\widehat{x},\widehat{y})|$$
	
	and $$|p_{1\varepsilon}+\theta_{2n}|\le  \frac{2}{\varepsilon}|\widehat{x}-\widehat{y}|+2K_n |\widehat{y}-y_n|+|\nabla_{\beta}^yP(\widehat{x},\widehat{y})|.$$
	Hence $p_{1\varepsilon}+\lambda q+\theta_{1n}, p_{1\varepsilon}+\theta_{2n}$ are contained in a bounded set. Let $R_1$ be an upper bound of $\widehat{x},\widehat{y}, p_{1\varepsilon}+\lambda q+\theta_{1n}$ and $p_{1\varepsilon}+\theta_{2n}.$

	Let $R\geq R_1,$ by the continuity at $0$ of the function $w_R: X^*_{\beta}\to \mathbb{R}$ in the assumption (H0), for every $\eta>0,$ there exists a neighborhood $V_{\beta}\subset X^*_\beta$ of 0 such that $|w_R(p)|<\eta$  for all $p\in V_{\beta}.$

	Since $\tau_F$  is the strongest topology  in the class of $\tau_{\beta}$ topologies on the space $X^*,$ there  exists $\delta>0$ such that $\delta B_1\subset \frac{1}{2}V_{\beta},$ where $B_1$ is the unit ball in $X^*.$

	On account of  both $K_n$ and $\delta_n$ converge to 0 as $n\to \infty,$ we can choose $n$ large enough such that $K_n.(\nabla_{\beta}|.|^2) (\widehat{x}-x_n)\in \delta B_1;$  
	$K_n.(\nabla_{\beta}|.|^2) (\widehat{y}-y_n)\in \delta B_1$ and $\delta_n<\delta.$ This implies that  $\theta_{1n}\in V_\beta$ and $ \theta_{2n}\in V_\beta.$

	Because $\Psi-P$ attains its maximum at $(\widehat{x},\widehat{y}),$ then $$\Psi(\widehat{x},\widehat{y}) \ge \Psi(x_n,y_n)+P(\widehat{x},\widehat{y})-P(x_n,y_n)=\sup_{\Delta(a)}\Phi-3\delta_n\ge 0$$
	for $n$ large enough.
	
	Combining with (\ref{1.16}),  (\ref{1.28}), we have
	\begin{equation}
	u(\widehat{x})-v(\widehat{y}) \ge  \Psi(\widehat{x},\widehat{y}) \ge \Psi(x_n,y_n)=\sup_{\Delta(a)}\Phi-3\delta_n\ge 0.
	\label{them7}
	\end{equation}
	Hence
	$u(\widehat{x})\ge v(\widehat{y}).$
	
	Since $H, \widehat{H}$ satisfy condition $(H0),$ for a sufficiently large $n,$ we have
	$$ |H(\widehat{x},u(\widehat{x}), p_{1\varepsilon}+\lambda q+\theta_{1n})-H(\widehat{x},u(\widehat{x}), p_{1\varepsilon}+\lambda q)|<|w (\theta_{1n}, R)|<\eta $$
	which implies
	$$ -H(\widehat{x},u(\widehat{x}), p_{1\varepsilon}+\lambda q+\theta_{1n})<-H(\widehat{x},u(\widehat{x}), p_{1\varepsilon}+\lambda q)+\eta  $$
	and
	$$ |\widehat{H}(\widehat{y},v(\widehat{y}),p_{1\varepsilon}+\theta_{2n})-\widehat{H}(\widehat{y},v(\widehat{y}),p_{1\varepsilon}) | <w (\theta_{1n}, R)<\eta $$
	It follows that
	$$ \widehat{H}(\widehat{y},v(\widehat{y}),p_{1\varepsilon}+\theta_{2n})< \widehat{H}(\widehat{y},v(\widehat{y}),p_{1\varepsilon})+\eta.$$
	This implies
	\begin{align}\label{bd1}
	\widehat{H}(\widehat{y},v(\widehat{y}),p_{1\varepsilon}+\theta_{2n})&-H(\widehat{x},u(\widehat{x}), p_{1\varepsilon}+\lambda q+\theta_{1n})\\
	&< \widehat{H}(\widehat{y},v(\widehat{y}),p_{1\varepsilon})-  H(\widehat{x},u(\widehat{x}), p_{1\varepsilon}+\lambda q)+2\eta.\notag
	\end{align}

	By (\ref{1.31'});  (\ref{1.31}) and (\ref{bd1}) we have
	\begin{equation*}
	u(\widehat{x})-v(\widehat{y})\le  \widehat{H}(\widehat{y},v(\widehat{y}),p_{1\varepsilon})-H(\widehat{x},u(\widehat{x}),p_{1\varepsilon}+\lambda q)+2\eta.
	\end{equation*}
	On the other hand, using (H1) (resp. (H1)*) we get 
	\begin{equation}\label{1.34}
	\begin{split} 
	\widehat{H}&(\widehat{y},v(\widehat{y}),p_{1\varepsilon})-H(\widehat{x},u(\widehat{x}),p_{1\varepsilon}+\lambda q)\\
	&= \widehat{H}(\widehat{y},v(\widehat{y}),p_{1\varepsilon})-H(\widehat{y},v(\widehat{y}),p_{1\varepsilon})+H(\widehat{y},v(\widehat{y}),p_{1\varepsilon})-H(\widehat{x},u(\widehat{x}),p_{1\varepsilon}+\lambda q) \\
	&\le  \widehat{H}(\widehat{y},v(\widehat{y}),p_{1\varepsilon})-H(\widehat{y},v(\widehat{y}),p_{1\varepsilon})+H(\widehat{y},u(\widehat{x}),p_{1\varepsilon})-H(\widehat{x},u(\widehat{x}),p_{1\varepsilon}+\lambda q)\\
	&=  \widehat{H}(\widehat{y},v(\widehat{y}),p_{1\varepsilon})-H(\widehat{y},v(\widehat{y}),p_{1\varepsilon})
	+H(\widehat{y},u(\widehat{x}),p_{1\varepsilon})-H(\widehat{x},u(\widehat{x}),p_{1\varepsilon})\\
	&+H(\widehat{x},u(\widehat{x}),p_{1\varepsilon})-H(\widehat{x},u(\widehat{x}),p_{1\varepsilon}+\lambda q)\\
	&\le\sup\{(\widehat{H}(z,r,p)-H(z,r,p))^+: (z,r,p)\in B_{R_0}\times \R\times X^*, |p|\le |p_{1\varepsilon}|\}\notag\\
	&+m_H(2|\widehat{x}-\widehat{y}|^2/\varepsilon+|\widehat{x}-\widehat{y}|)+\sigma_H(\lambda,2|\widehat{x}-\widehat{y}|/\varepsilon+\lambda);\\
	\Big({\rm resp.}&\quad  \widehat{H}(\widehat{y},v(\widehat{y}),p_{1\varepsilon})-H(\widehat{x},u(\widehat{x}),p_{1\varepsilon}+\lambda q)\\
	&\leq \sup\{(\widehat{H}(z,r,p)-H(z,r,p))^+: (z,r,p)\in B_{R_0}\times \R\times X^*, |p|\le |p_{1\varepsilon}|\}\\
	&+m_H(2|\widehat{x}-\widehat{y}|^2/\varepsilon+|\widehat{x}-\widehat{y}|)+\sigma_H(\lambda,2|\widehat{x}-\widehat{y}|/\varepsilon+\lambda)+L_H(u(\widehat{x})-v(\widehat{y}))	\Big),	
	\end{split}
	\end{equation}
	where $B_{R_0}$  is the ball centered at 0 with  sufficiently large radius in $X.$  From \eqref{1.24} and $(x_n,y_n)\in S_n$
	we have
	\begin{equation}
	|\widehat{x}-\widehat{y}|\le |x_n-y_n|+2\gamma\le (\varepsilon m(a))^{1/2}+2\gamma.
	\label{1.36}
	\end{equation}
	It follows that
	\begin{align*}
	2|\widehat{x}-\widehat{y}|^2/\varepsilon+|\widehat{x}-\widehat{y}|&\le \frac{2}{\varepsilon}(\varepsilon m(a)+4\gamma (\varepsilon m(a))^{1/2} +4\gamma^2)+(\varepsilon m(a))^{1/2}+2\gamma\\ 
	&=2 m(a)+(\varepsilon m(a))^{1/2}+\gamma\left(\frac{8}{\varepsilon}(\varepsilon m(a))^{1/2}+\frac{8\gamma^2}{\varepsilon} +2\right).
	\end{align*}
	Taking $\gamma\to 0,$ we arrive at  
	\begin{equation*}
	2|\widehat{x}-\widehat{y}|^2/\varepsilon+|\widehat{x}-\widehat{y}|\le 2 m(a)+(\varepsilon m(a))^{1/2}.
	\label{them6}
	\end{equation*}
	From (\ref{1.32}) and (\ref{1.36})  we have
	\begin{equation}
	|p_{i\varepsilon}|\le 2( (\varepsilon m(a))^{1/2}+\gamma)/\varepsilon.
	\label{1.37}
	\end{equation}
	From  (\ref{them7})-(\ref{1.37}) we have
	\[\Phi(x_n,y_n) -3\delta_n \le  C(a,\lambda,\varepsilon)+2\eta\]
	\[\Big({\rm resp.}\quad  \Phi(x_n,y_n) -3\delta_n \le \frac{1}{1-L_H}( C(a,\lambda,\varepsilon)+2\eta) \Big).\]
	
	Let $n\to \infty$ and $\eta\to0$ we have \eqref{them1}, and the inequality \eqref{1.19} is proved completely.
	
	Finally, choose $R\geq \max\{R_0,R_1\},$ where $R_0$ and $R_1$ appear in Step 2. By the definition of  $\zeta,$ for $|x|<R,$ $\zeta(|x|-R)=0$ and
	\begin{equation}\label{4.41}
	u(x)-v(x)=\Phi(x,x).
	\end{equation}
	Because $(x,x)\in\Delta(a),$ so from \eqref{1.19} we obtain 
	\begin{equation}\label{them9}
	\begin{split}
	&u(x)-v(x)\le \sup_{\partial \Omega} (u-v)^++2m(a)+C(a,\lambda,\varepsilon)\\
	\Big({\rm resp.}&\quad u(x)-v(x)\le \sup_{\partial \Omega} (u-v)^++2m(a)+\frac{1}{1-L_H}C(a,\lambda,\varepsilon)  \Big).
	\end{split}
	\end{equation}
	
	Thank to \eqref{4.41},  $u-v$ is a bounded function. Thus, we can choose $B=0$ in \eqref{1.14}. By taking the limit in \eqref{them9} as $\lambda\to 0$ and $a\to 0$, we obtain \eqref{1.11}.

	\subsection{The case $\Omega=X$}
	In this case we have $\partial \Omega=\emptyset,$  so that $\sup_{\partial \Omega} (u-v)^+=0,$  $\rho(x)=\infty$ for all $x\in X.$  Using the same arguments as above, we arrive at the conclusion of the theorem. The theorem is completely proved.

\bigskip
	
\noindent {\bf Acknowledgements}

The authors thank an anonymous referee for careful reading
of the manuscript and her/his useful comments leading to an
improvement of the presentation.

\end{document}